\documentclass{amsart}	
\usepackage[utf8]{inputenc}

\usepackage{fullpage}

\usepackage{array}

\usepackage{amssymb, amsmath, amsthm, color, amsfonts, caption, url, longtable, mathtools, hyperref}
\usepackage{epsfig}

\usepackage{comment}
\usepackage{hyperref}
\usepackage{url}
\newcommand{\bburl}[1]{\textcolor{blue}{\url{#1}}}

\newtheorem{theorem}{Theorem}
\numberwithin{theorem}{section}

\newtheorem{lemma}[theorem]{Lemma}

\newtheorem{example}[theorem]{Example}

\newcommand{\N}{\mathbb{N}}

\newcommand{\R}{\mathbb{R}}

\title{Sums of Powers by L'Hopital's Rule}
\author{Eduardo Dueñez, Asimina S. Hamakiotes and Steven J. Miller}

\address{Department of Mathematics, University of Texas at San Antonio, TX 78249\\ \ } \email{eduardo.duenez@utsa.edu\\ \ }

\address{Department of Mathematics, University of Connecticut, Storrs, CT 06269\\ \ } \email{asimina.hamakiotes@uconn.edu\\ \ }

\address{Department of Mathematics and Statistics, Williams College, Williamstown, MA 01267} \email{sjm1@williams.edu}

\thanks{We thank our students and colleague \'Alvaro Lozano-Robledo for helpful discussions.}



\date{\today}

\begin{document}

\maketitle


\vspace{-0.75cm}

\begin{abstract}
For a positive integer $d$, let $p_d(n) := 0^d + 1^d + 2^d + \cdots + n^d$; i.e., $p_d(n)$ is the sum of the first $d$\textsuperscript{th}-powers up to $n$. It's well known that $p_d(n)$ is a polynomial of degree $d+1$ in $n$. While this is usually proved by induction, once $d$ is not small it's a challenge as one needs to know the polynomial for the inductive step. We show how this difficulty can be bypassed by giving a simple proof that $p_d(n)$ is a polynomial of degree $d+1$ in $n$ by using L'Hopital's rule, and show how we can then determine the coefficients by Cramer's rule. This illustrates a general principle and the point of our paper: there's more than one path to a goal, different approaches have their advantages and disadvantages, and the more techniques one knows, the more likely one can successfully attack a problem.
\end{abstract}




\section{Introduction}

One of the first problems that arise when we encounter a sequence of numbers is to add them, and frequently beautiful, simple formulas emerge. For example, the geometric series with ratio $r$ and zeroth term $a_0$ has $n$\textsuperscript{th} term $a_n = a_0 r^n$, and thus $$a_0 + a_1 + \cdots + a_n \ = \ a_0(1 + r + r^2 + \cdots + r^n) \ =: \ a_0 S_r(n).$$ Thus we just need to determine the answer in the simpler case when $a_0 = 1$. The standard trick is multiply by $r$ and subtract, and note we have a telescoping series: \begin{eqnarray} S_r(n) & \ = \ & 1 + r + r^2 + \cdots + r^n \nonumber\\ \underline{\ \ \ \ \ \ \ r S_r(n)} & \ = \ & \underline{\textcolor{white}{1 + }\ r + r^2 + \cdots + r^n + r^{n+1}} \nonumber\\ (1-r) S_r(n) & \ = \ & 1 \textcolor{white}{+ r + r^2 + \cdots + r^n\ } - r^{n+1}. \nonumber \end{eqnarray} Thus $$S_r(n) \ = \ 1 + r + \cdots + r^n \ = \ \frac{1 - r^{n+1}}{1-r} \ = \ \frac{1}{1-r} - \frac{r^{n+1}}{1-r},$$ and we have a formula for the sum of a finite geometric series. If $|r| < 1$ we can take the limit as $n \to \infty$ and find that the infinite geometric series has sum $1/(1-r)$. As the primary goal of this article is to highlight alternative approaches, we refer the reader to \cite{M} for another approach using a basketball game and memoryless processes.\footnote{See also the video from the third named author's first lecture in his probability courses: \bburl{https://youtu.be/aD7tIO0CGyg} (35:19); slides available here: \bburl{https://web.williams.edu/Mathematics/sjmiller/public_html/341Fa21/talks/firstlecture341.pdf}.}

\begin{figure}
\begin{center}
\scalebox{1}{\includegraphics{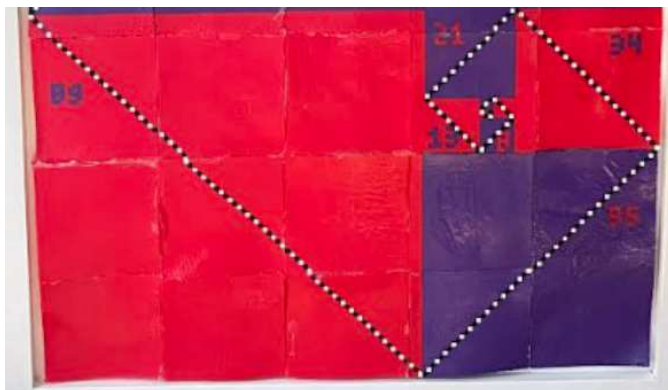}}
\caption{\label{fig:FibonacciSpiral55} A fuse bead tiling (constructed by Kayla and Steven Miller) of the plane by squares whose lengths are the Fibonacci numbers.}
\end{center}\end{figure}

Once we know the geometric series formula, many other sums can be done (including the main sums in this paper: power sums). For example, defining the Fibonacci numbers by $F_0 = 0, F_1 = 1$ and $F_{n+1} = F_n + F_{n-1}$ one can prove Binet's Formula\footnote{Not surprisingly there are many proofs of this result, from induction to generating functions as well as interesting combinatorial proofs (see \cite{BQ}).}: $$F_n \ = \ \frac1{\sqrt{5}} \left(\frac{1+\sqrt{5}}{2}\right)^n - \frac1{\sqrt{5}} \left(\frac{1-\sqrt{5}}{2}\right)^n.$$ By repeated applications we can determine sums of any power of Fibonacci numbers, though this approach may lead to an expression that is not as simple as possible; thus there may be a better avenue to these results! For example, $F_1^2 + F_2^2 + \cdots + F_n^2 = F_n \cdot F_{n+1}$. We can prove this using the geometric series formula and squaring, we can also prove it by induction, but interestingly there's an old geometric proof which the third named author has used numerous times in math outreach activities with kids as young as kindergarten. For each $n \le 13$, groups of students are given exactly one $n \times n$ square, and told to place as many possible on the table as they can, with the rules that all must be flat on the table, and as they add squares one at a time at every moment the shape must be a rectangle. The groups quickly realize you cannot have more than one square as the side lengths don't line up, so to have a chance they need more. The smallest addition we can give is one extra $1\times 1$ square, and then every group discovers the Fibonacci sequence that covers the plane; see Figure \ref{fig:FibonacciSpiral55}.

We can compute the area of the rectangle two ways. First, it's just length times width, or $89 \cdot 144 = F_{11} \cdot F_{12}$. Second, it's the sum of the areas of the squares, or $F_1^2 + F_2^2 + \cdots + F_{12}^2$, proving our formula!

The next common sums are arithmetic progressions, where each term is $d$ more than the previous. Thus if we again start with $a_0$ as the zeroth term, now $a_n = a_0 + nd$ we find after some simple algebra $$a_0 + a_1 + \cdots + a_n\ = \ (n+1)a_0 + d(0 + 1 + 2 + \cdots + n).$$ Thus determining its sum is equivalent to summing the integers up to $n$.

More generally, we consider the problem of computing the $d$\textsuperscript{th}-power sums for $d$ a positive integer; this is often encountered in classes as a way to test the Fundamental Theorem of Calculus by evaluating the limit of the Riemann sums of $f(x) = x^d$. We define $p_d(n) \coloneqq 0^d + 1^d + 2^d + \cdots + n^d$ to be the sum of the first $d$\textsuperscript{th}-powers up to $n$. Many students are introduced to these problems through a story of a young Gauss, whose teacher once had a bad day and asked the class to compute $p_1(100) = 0 + 1 + 2 + \cdots + 100$, expecting to have some time to himself while the students tediously added term after term. Gauss however quickly noticed that if the sum is written in reverse order and added to itself, the result is
\begin{align*}
    2p_1(100) \ = \ (100 + 0) + (99+1) + (98+2) + \cdots + (0+100) \ = \ 100\cdot 101,
\end{align*}
which leads to $p_1(100) = 50\cdot 101 = 5050$. More generally, this beautiful trick yields
\[
p_1(n) \ := \ 0 + 1 + 2 + \cdots + n \ = \ \frac{n(n+1)}{2} \ = \ \frac{1}{2}n^2+ \frac{1}{2}n,
\]
and we can notice a lot of interesting items from this formula.

\begin{enumerate}

    \item First, $p_1(n)$ is a polynomial of degree one higher than the power we are summing; i.e., $p_1(n)$ sums the first powers and is a polynomial of degree two.

    \item Second, from $p_1(n) = (1/2)n^2+ (1/2)n$, we notice that the coefficients are not integers. The leading term has coefficient $1/2$ and there's no constant term.
\end{enumerate}

It's easy to see both of these hold for additional $d$. Using a straightforward inductive argument we can verify that for the sum of squares
\begin{align*}
p_2(n) \ = \ 0^2 + 1^2 + 2^2 + \cdots + n^2 \ = \ \frac{n(n+1)(2n+1)}{6} \ = \ \frac{1}{3}n^3 + \frac{1}{2}n^2 + \frac{1}{6}n,
\end{align*}
and the sum of third powers is
\begin{align*}
p_3(n) \ = \ 0^3 + 1^3 + 2^3 + \cdots + n^3 \ = \  \left(\frac{n(n+1)}{2}\right)^2 \ = \ \frac{1}{4}n^4 + \frac{1}{2}n^3 + \frac{1}{4}n^2.
\end{align*}
Observe that $p_2(n)$ and $p_3(n)$ also have no constant term, and the coefficient of the leading term is $1/(d+1)$ (the leading coefficient of $p_2(n)$ is $1/3$ and the leading coefficient of $p_3(n)$ is $1/4$). Given these observations, it's natural to conjecture that $p_d(n)$ is a polynomial of degree $d+1$ with no constant term and leading term $1/(d+1)$. The reason these are easy to prove by induction for small $d$ is because we have a conjectured formula for the polynomial. If we did not know the polynomial, we would first have to guess it. This can be done by assuming a polynomial relation, evaluating for $n = 0, 1, 2, \dots, d$ and then interpolating to see what polynomial; see \S\ref{matrixmethod}.

Note that if $p_d(n)$ is a polynomial then the main term must be $n^{d+1}/(d+1)$; we can see this by approximating the sum with an integral:
\begin{align*}
    0^d + 1^d + \cdots + n^d \ \approx \ \int_0^n x^ddx \ = \  \frac{x^{d+1}}{d+1}\Bigg\vert_0^n \ = \ \frac{1}{d+1}n^{d+1}.
\end{align*}
With a little more effort, we can obtain upper and lower bounds for how close the approximation is.

With even more work, we can use the Euler-Maclaurin formula (see for example Section 9.5 in \cite{GKP}) to prove that $p_d(n)$ is a polynomial. For a function $f(x)$ that is $k$ times continuously differentiable on the interval $[0,n]$, the difference between the sum and the integral of $f(x)$ is
\[
\sum_{m=0}^n f(m) - \int_0^n f(x)dx \ = \  \sum_{j=1}^k \frac{B_j}{j!}\left(f^{(j-1)}(n) - f^{(j-1)}(0) \right) + R_k,
\]
where $B_j$ is the $j$\textsuperscript{th} Bernoulli number and $R_k$ is the error term. The error term is in terms of periodized Bernoulli functions, $P_j(x) = B_j(x - \lfloor x\rfloor)$, and is
\[
R_k \ =\ (-1)^{k+1}\int_0^n f^{(k)}(x) \frac{P_k(x)}{k!}dx.
\]
If we take $d+1$ derivatives, then there's no error term as $f(x) = x^d$ and thus $p_d(n)$ is a polynomial in $n$. We can derive the Euler-Maclaurin formula with the error term by repeatedly applying integration by parts on
\begin{align*}
R_k \ = \ (-1)^{k+1}\int_0^n f^{(k)}(x) \frac{P_k(x)}{k!}dx \ = \ (-1)^{k+1} \left( f^{(k+1)}(n) \frac{P_k(n)}{k!} - f^{(k+1)}(0) \frac{P_k(0)}{k!}\right) + R_{k-1},
\end{align*}
for $k\geq 2$. The boundary terms contribute to the main terms in the Euler-Maclaurin formula and the remainder integrals contribute to the error term. Therefore, using $f(x) = x^d$, we quickly obtain an integration proof that $p_d(n)$ is a polynomial.


Our goal is to give a new, simple proof (which we have not been able to find in the literature) that $p_d(n)$ must be a polynomial by using L'Hopital's rule, bypassing the need to use an inductive argument. Thus the differentiation argument below complements the integration approach above, and highlights a key fact of mathematics: often there's more than one way to an answer, and different approaches have their advantages and disadvantages.

\begin{theorem}\label{MainThm}
Let $p_d(n)$ be the sum of the $d$\textsuperscript{th}-powers up to $n$. Then L'Hopital's rule implies that $p_d(n)$ is a polynomial of degree $d+1$ with leading term $1/(d+1)$ and constant term zero.
\end{theorem}

We can find this polynomial using the following procedure; we give the details in the next section. Start with the finite geometric series formula
\begin{align*}
    1 + x + x^2 + \cdots + x^n\ =\ \frac{x^{n+1}-1}{x-1},
\end{align*}
and apply $x\frac{d}{dx}$ to both sides $d$ times, then evaluate the left hand side at $x=1$ and take the limit as $x$ approaches $1$ of the right hand side. This involves applying L'Hopital's rule $2^d$ times. While we can explicitly determine the coefficients, as doing so becomes computationally involved for even modest $d$ we discuss in \S\ref{matrixmethod} how to find the coefficients, knowing it's a polynomial, through Cramer's rule. Using a similar method but applying $\frac{d}{dx}$ instead, we obtain similar results for related sums, which are interesting as well.


\section{Differentiation Method for Computing $p_d(n)$}\label{polynomial}

Our method yields an easy, new way to prove $p_d(n)$ has to be a polynomial, and can be used to explicitly find that polynomial. Unfortunately in practice it becomes computationally intensive once $d$ is large; $d=2$ is already a significant amount of work, as we see below.

As remarked we start with the geometric series formula for a finite sum, proved in the introduction; note there are no conditions on $x$ other than it doesn't equal 1 (in that case, the sum on the right hand side is just $n+1$).

\begin{lemma}[Finite Geometric Series Formula]
For $x \neq 1$, we have
\begin{align*}
     1 + x + x^2 + \cdots + x^n \ = \ \frac{x^{n+1}-1}{x-1}.
\end{align*}
\end{lemma}

Since $S_n(x)$ is a finite sum, we do not need to worry about whether or not $|x|<1$, though we do need $x\neq 1$ so that the denominator is not zero; of course if $|x| < 1$ we can take the limit and obtain the formula for the infinite geometric sum.

Since the sum is finite, we can differentiate term by term using the derivative of a finite sum is the sum of the derivatives. This can easily be proved by induction; it's not always true for infinite sums. We start with the geometric series
\[
1 + x + x^2 + \cdots + x^n \ = \ \frac{x^{n+1}-1}{x-1},
\]
and apply the operator $x\frac{d}{dx}$; we could apply just $\frac{d}{dx}$, but our choice keeps the exponent of each term unchanged, and this will be useful as we take repeated derivatives in proving the result for general $d$. We find
\begin{align*}
0 + 1x + 2x^2 + \cdots + nx^n \ &= \ x \left[ \frac{(n+1)x^n(x-1)-(x^{n+1}-1)}{(x-1)^2}\right], \\
&= \ \frac{nx^{n+2}-(n+1)x^{n+1}+x}{(x-1)^2},
\end{align*}
where the quotient rule is used to differentiate the right hand side. Note the denominator is of degree 2 in $x-1$.%
\footnote{\label{fn:1}Since
\begin{equation*}
  x\frac{d}{dx}\frac{f(x)}{(x-1)^k}
  = \frac{x(x-1)f'(x) +kxf(x)}{(x-1)^{k+1}}
\end{equation*}
for any polynomial $f(x)$, it's clear that $d$ successive applications of the operator $x(d/dx)$ to the geometric series yields an expression of the form $g(x)/(x-1)^{d+1}$.}

Observe that if we take $x=1$, the left hand side becomes the sum of first powers up to $n$
\[
0 + 1(1) + 2(1)^2 + \cdots + n(1)^n \ = \ 0 + 1 + 2 + \cdots + n,
\]
while the right hand side becomes $0/0$:
\[
\frac{n(1)^{n+2}-(n+1)(1)^{n+1}+(1)}{((1)-1)^2} \ = \ \frac{n-n-1+1}{0} \ = \ \frac{0}{0}.
\]
We use L'Hopital's rule to evaluate the right hand side. We know that we must obtain something nice on the right hand side as the left hand side is nice. If we take the limit as $x$ approaches 1 and apply L'Hopital twice, we get
\begin{align*}
    \lim_{x\to 1} \frac{nx^{n+2}-(n+1)x^{n+1}+x}{(x-1)^2} \ &= \ \lim_{x\to 1} \frac{n(n+2)x^{n+1}-(n+1)(n+1)x^{n}+1}{2(x-1)}, \\
    &= \ \lim_{x\to 1} \frac{n(n+2)(n+1)x^{n}-(n+1)^2nx^{n-1}}{2}, \\
    &= \ \frac{(n^3+3n^2+2n)-(n^3+2n^2+n)}{2}, \\
    &= \ \frac{n^2 + n}{2},
\end{align*}
which is the claimed formula $p_1(n) = n(n+1)/2$.

Before we apply L'Hopital's rule, we can try to figure out what the right hand side should be. As remarked earlier we know from the Riemann sum
\[
0 + 1 + 2 + \cdots + n \ \approx \ \int_0^n xdx \ =\ \frac{1}{2}n^2
\]
that the right hand side has to start with leading term $n^2/2$. By taking $n=0$ we see that if it's a polynomial then the constant term must be zero. We also know that there cannot be any $n$-dependence from the denominator, because the denominator is always going to be powers of $(x-1)$ and as we apply $x\frac{d}{dx}$ and L'Hopital's rule we never have any $n$-dependence downstairs. As for the numerator, once we take sufficiently many derivatives so that the denominator is a constant we see that when we take $x=1$ it doesn't matter what the exponents of $x$ are, because we are just going to get 1. Thus upstairs are polynomials in $n$. Therefore, after two applications of L'Hopital's rule we get that the right hand side is a polynomial in $n$ of degree 2.

We now find the formula for the sum of squares, though the algebra is less pleasant. As before, we first apply  $x\frac{d}{dx}$ to both sides of the finite geometric series formula:
\begin{align*}
0 + 1x + 2x^2 + \cdots + nx^n \ &= \  \frac{nx^{n+2}-(n+1)x^{n+1}+x}{(x-1)^2}, \end{align*} and then applying $x\frac{d}{dx}$ again yields
\begin{align*}
x + 2^2x^2 + 3^2x^3 + \cdots + n^2x^n \ &= \
\frac{n^2x^{n+4}+(-3n^2-2n+1)x^{n+3}+(3n^2+4n)x^{n+2}-(n+1)^2x^{n+1}-x^3 +x}{(x-1)^4}.
\end{align*}
Taking $x=1$ yields the left hand side is the sum of squares up to $n$, whereas the right hand side becomes $0/0$:
\begin{multline*}
\frac{n^2(1)^{n+4}+(-3n^2-2n+1)(1)^{n+3}+(3n^2+4n)(1)^{n+2}-(n+1)^2(1)^{n+1}-(1)^3 +1}{(1-1)^4} \\  = \  \frac{n^2+(-3n^2-2n+1)+(3n^2+4n)-(n+1)^2-1 +1}{(1-1)^4} \ = \ \frac{0}{0}.
\end{multline*}
Observe that the denominator of the right hand side is $(x-1)^4$, which has no $n$-dependence. The only $n$-dependence is in the numerator. The right hand side is a polynomial after we take the limit because every time we apply L'Hopital's rule, we obtain polynomials of $n$ coming down from the exponents of the powers of $x$. When we evaluate at $x=1$, the monomials $x^k$ contribute 1 no matter what their exponents are, and we are left with a polynomial in $n$ of degree 3. We have to be careful about how things cancel, but since we know that the left hand side is nice, we know after some algebra that the right hand side must have a nice limit as well. Further, as we know from the Riemann sum calculation what the leading term must be, we know it'll be a polynomial of degree $d+1$.

After taking the limit as $x$ approaches 1 of the right hand side and applying L'Hopital four times, we get
\begin{align*}
    \lim_{x\to 1}\frac{n^2x^{n+4}+(-3n^2-2n+1)x^{n+3}+(3n^2+4n)x^{n+2}-(n+1)^2x^{n+1}-x^3 +x}{(x-1)^4} \ = \ \frac{2n^3+3n^2+n}{6},
\end{align*}
which is the desired $p_2(n) = n(n+1)(2n+1)/6$. Note we know we must apply L'Hopital's rule exactly four times, as the limit exists and fewer times results in a zero denominator.

We can keep applying $x\frac{d}{dx}$ and L'Hopital's rule to get formulas for $p_3(n), p_4(n), p_5(n)$ and so on. Not only will we find that these sums are polynomials, but we'll obtain their coefficients. Unfortunately we need to use L'Hopital $2^d$ times for the sum of the $d$\textsuperscript{th} powers, and the resulting algebra, though doable, becomes unpleasant. It's thus better to simply note that this argument immediately yields that it must be some polynomial with zero constant term and leading term $x^{d+1}/(d+1)$, and find another way to determine the coefficients. This can be done by Cramer's rule; the argument is standard, but given in \S\ref{matrixmethod} for completeness. We next turn to analyzing what happens if we apply $\frac{d}{dx}$ instead of $x \frac{d}{dx}$.


\section{Algebraic remarks}\label{sec:alg-rem}

Recall the ``falling powers'' defined by
\begin{equation*}
  (y)_k \ := \  \prod_{i=1}^k(y-i+1) \ = \  y (y-1) \cdots (y-k+1)
\end{equation*}
for $k\in\N$ (with the usual ``empty product'' convention $(y)_0=1$).
Let
$$f_n(x)\ \coloneqq\ 1+x+x^2+\cdots+x^n\ = \ \frac{x^{n+1}-1}{x-1},$$ and
\begin{align*}
  g_k(x) &\ \coloneqq\ \left( x\frac{d}{dx} \right)^kf_n(x)
           \ = \  x + 2^kx^2 + 3^kx^3 + \cdots + n^kx^n, \\
  h_k(x) &\ \coloneqq\ \frac{d^kf_{n+k}}{dx^k}
   \ = \  (k)_k + (k+1)_kx + (k+2)_kx^2 + (k+3)_kx^3 + \cdots + (n+k)_kx^n.
\end{align*}
Note that $g_0(x) = 1 + x + \cdots + x^n$ starts with the term $1 = 0^0$; as written, the right-hand expression for $g_0$ above holds only for $k\ge 1$ (adding the funny-looking term $0^k$ to that expression makes it always correct).
The sum of the $(n+1)$ first ``falling $k$-th powers'',
\begin{equation*}
  q_k(n) \ \coloneqq \ h_k(1) \ = \  (k)_k + (k+1)_k + \cdots + (n+k)_k,
\end{equation*}
has the simple closed form $q_k(n) = (n+k+1)_k/(k+1)$. This can easily be proved by induction.
In a manner similar to our earlier arguments, we can also use repeated applications of L'Hopital's rule to show it's a polynomial of degree $k+1$ in $n$.

For fixed $k$, a moment's reflection shows that $\{g_0, \dots, g_k\}$ and $\{h_0, \dots, h_k\}$ span (over $\R$, say) the same linear space of polynomials.
In fact, this follows from the even simpler observation that the falling powers $(y)_0, (y)_1, \dots, (y)_k$ of an indeterminate $y$ have the same span---consisting of polynomials of degree at most $k$---as the monomials $1, y, y^2, \dots, y^k$.
Thus, implicit in the differentiation method of Section~\ref{polynomial} is a method for performing the change of basis.
Interested readers are invited to make the construction explicit, expressing the standard monomials $y^k$ in the form
\begin{equation*}
  y^k \ = \  \sum_{j=0}^k S(k,j) (y)_k,
\end{equation*}
thus yielding a (well known otherwise) recursive definition for the coefficients $S(k,j)$ (which are Stirling numbers of the second kind).%
\footnote{The standard combinatorial definition of $S(k,j)$ is the number of partitions of an $n$-element set into $k$ (nonempty) parts.}


\section{Matrix Method for Computing $p_d(n)$}\label{matrixmethod}

By Theorem \ref{MainThm}, we know that $p_d(n)$ is a polynomial of degree $d+1$ with leading term $n^{d+1}/(d+1)$ and constant term zero. As remarked, the standard way to prove this is by induction, see for example \cite{MT-B}, but the difficulty with induction is that we need to know what our target is. Specifically, we need to have a polynomial for $p_d(n)$. Our L'Hopital approach immediately proved the sum is a polynomial, but determining the coefficients is a tedious calculation. We show how to find the coefficients easily by interpolation; note this could also be used to set up the proof by induction.



Suppose that $p_d(n) = a_1n + \cdots + a_dn^d + a_{d+1}n^{d+1}$. Theorem \ref{MainThm} tells us that $a_{d+1} = 1/(d+1)$ and $a_0=0$. If we evaluate $p_d(n)$ at $n = 1,2,\ldots, d+1$, then we get the following $d+1$ system of equations:
\begin{align*}
    p_d(1) \ &= \ a_1 1 + \cdots +  a_d 1^d + a_{d+1} 1^{d+1}\ = \ 0^d + 1^d, \\
    p_d(2) \ &= \ a_1 2 + \cdots + a_d 2^d + a_{d+1} 2^{d+1} \ = \ 0^d + 1^d + 2^d, \\
    &\ \vdots \\
    p_d(d+1) \ &= \ a_1(d+1) + \cdots + a_d(d+1)^d + a_{d+1}(d+1)^{d+1} \ = \ 0^d + 1^d + \cdots + (d+1)^d.
\end{align*}
We can set up the following in matrix form in order to solve for the coefficients
\begin{align*}
    \begin{pmatrix}
    1&1&\cdots & 1\\
    2&2^2&\cdots& 2^{d+1}\\
    \vdots&\vdots&\cdots&\vdots\\
    (d+1)&(d+1)^2&\cdots&(d+1)^{d+1}
    \end{pmatrix}
    \begin{pmatrix}
    a_1\\a_2\\\vdots\\a_{d+1}
    \end{pmatrix} \ = \
    \begin{pmatrix}
    p_d(1) \ = \ 0^d + 1^d\\p_d(2) \ = \ 0^d + 1^d + 2^d\\\vdots\\p_d(d+1) \ = \ 0^d + 1^d + \cdots + (d+1)^d
    \end{pmatrix}.
\end{align*}
Observe that the square matrix on the left hand side is a Vandermonde matrix and is invertible as all of the columns are distinct; for more on the Vandermonde matrix, see for example Section 2.8.1 in \cite{PTVF}. We can solve for the $a_i$'s using Cramer's rule as follows:
\begin{align*}
    \begin{pmatrix}
    a_1\\a_2\\\vdots\\a_{d+1}
    \end{pmatrix} \ = \
    \begin{pmatrix}
    1&1&\cdots & 1\\
    2&2^2&\cdots& 2^{d+1}\\
    \vdots&\vdots&\cdots&\vdots\\
    (d+1)&(d+1)^2&\cdots&(d+1)^{d+1}
    \end{pmatrix}^{-1}
    \begin{pmatrix}
    p_d(1)\\p_d(2)\\\vdots\\p_d(d+1)
    \end{pmatrix}.
\end{align*}
Now that we have the coefficients, we can substitute them into our formula $p_d(n) = a_{d+1}n^{d+1}+a_dn^d + \cdots + a_1n$ and find the polynomial. Thus, if $p_d(n)$ is a polynomial, then there's only one polynomial that it could be, and we have a way to find it (though it's not pleasant for even modest $n$).


We conclude this section with an example of using this method.

\begin{example}[Using matrix method to find the formula for $p_1(n)$]
Recall that $p_1(n) = 0 + 1 + 2 + \cdots + n$. Let $p_1(n) = a_2n^2 + a_1n$, as we know trivially there's no constant term. By setting up the following
\begin{align*}
    \begin{pmatrix}
    1&1\\
    2&4\\
    \end{pmatrix}
    \begin{pmatrix}
    a_1\\a_2
    \end{pmatrix} \ = \
    \begin{pmatrix}
    p_1(1) = 0 + 1\\p_1(2) = 0 + 1 + 2
    \end{pmatrix} \ = \
    \begin{pmatrix}
    1\\3
    \end{pmatrix},
\end{align*}
we can solve for the coefficients
\begin{align*}
    \begin{pmatrix}
    a_1\\a_2
    \end{pmatrix} \ = \
    \begin{pmatrix}
    1&1\\
    2&4\\
    \end{pmatrix}^{-1}
    \begin{pmatrix}
    1\\3
    \end{pmatrix} \ = \
    \frac{1}{2}\begin{pmatrix}
    4&-1\\
    -2&1\\
    \end{pmatrix}
    \begin{pmatrix}
    1\\3
    \end{pmatrix} \ = \
    \begin{pmatrix}
    1/2\\1/2
    \end{pmatrix}.
\end{align*}
We find $a_2 = 1/2$ and $a_1 = 1/2$, yielding $p_1(n) = (1/2)n^2 + (1/2)n = n(n+1)/2$.
\end{example}


\section{Related Sums}


We saw what happens when we apply $x\frac{d}{dx}$ to the finite geometric series: we get sums of powers. If we apply $\frac{d}{dx}$ instead, we obtain related sums; it's a standard problem in combinatorics to pass from one sum to a related one. Frequently the challenge is to determine the right trade-off: is it better to do some combinatorics or a harder computation?



We again start with the geometric series
\begin{align*}
1+x+x^2+\cdots+x^n \ = \ \frac{x^{n+1}-1}{x-1}.
\end{align*}
If we apply $\frac{d}{dx}$ to both sides, we get
\begin{align*}
    0 + 1 + 2x + 3x^2 + \cdots + nx^{n-1} \ = \ \frac{nx^{n+1}-(n+1)x^n + 1}{(x-1)^2},
\end{align*} note as we are not applying $x \frac{d}{dx}$ the degree of the numerator is smaller, and thus the derivatives are a little simpler.

If we apply $\frac{d}{dx}$ to both sides again, we get
\begin{multline*}
    0 + 2 + 6x + \cdots + n(n-1)x^{n-2} \\ = \ \frac{(x-1)^2(n(n+1)x^n -n(n+1)x^{n-1})-(nx^{n+1}-(n+1)x^n+1)(2(x-1))}{(x-1)^4}.
\end{multline*}
Unfortunately, when we take $x=1$, we can see from evaluating the left hand side, that we no longer have the sum of squares. Instead, we get
\[
1(0) + 2(1) + 3(2) + \cdots + n(n-1) \ = \ \sum_{k=1}^n k(k-1),
\]
where we can get a formula for this sum by taking the limit as $x$ approaches 1 of the right hand side and applying L'Hopital like we did before.

Although this doesn't result in the sum of squares, we can use this to solve for the sum of squares. Observe that we can rewrite the summation above as follows,
\[
\sum_{k=1}^n k(k-1) \ = \ \sum_{k=1}^n k^2 - \sum_{l=1}^{n} l.
\]
Since we know $\sum_{l=1}^{n}l$ is $p_1(n)$, the sum of first powers up to $n$, we can determine $\sum_{k=1}^n k^2$ is $p_2(n)$, the sum of squares up to $n$. If we know that $\sum_{l=1}^{n}l$ is a polynomial, which we do, then we know that $\sum_{k=1}^n k^2$ has to be a polynomial as well, and we can solve for that polynomial since we have polynomials for $\sum_{k=1}^n k(k-1)$ and $\sum_{l=1}^{n}l$. We can generalize this approach and similarly show by induction that $p_d(n)$ is a polynomial, and find a formula for its coefficients.



\section{Concluding Remarks}

The goal of this note is to highlight that there \emph{can be} more than one way to prove a result. For some problems, frequently one technique works better than another, so the more ideas you see, the more likely you are to have something in your arsenal that works on future problems. Our starting point is the geometric series formula; while we gave the standard telescoping sum proof, as remarked see \cite{M} and the video mentioned in the footnote for an alternative approach using game theory and memoryless processes, further supporting the point of this paper. In general it's hard to prove identities; the methods here show how to generate infinitely many more, which is tremendously valuable, using just simple differentiation. Differentiating identities is a very important technique in probability, and is frequently used to determine the moments of a probability distribution.\footnote{The $k$\textsuperscript{th} moment of a continuous density $p(x)$, which is a non-negative function which integrates to 1 (there is of course a similar formulation for discrete densities), is $\int_{-\infty}^\infty x^k p(x) dx$. Note for some distributions these exist for all $k$, such as the uniform or Gaussian, while for others they exist for only finitely many, such as the Cauchy.} For more on this, see \cite{M} and the corresponding online lectures from the third named author's probability class based on this book, for example \bburl{https://web.williams.edu/Mathematics/sjmiller/public_html/341Fa21/}, in particular Lectures 18 and 19.


\ \\



\end{document}